\documentstyle[psamsfonts,amssymb,amscd,12pt]{amsart}

\newtheorem{defn}{Definition}
\newtheorem{thm}[defn]{Theorem}

\newtheorem{cor}[defn]{Corollary}
\newtheorem{lem}[defn]{Lemma}
\newtheorem{prop}[defn]{Proposition}

\theoremstyle{remark}
\newtheorem{rem}{Remark}
\theoremstyle{remark}
\newtheorem{exam}[rem]{Example}

\numberwithin{equation}{section}
\numberwithin{defn}{section}

\begin{document}


\newcommand\spk{{\operatorname{Spec}}(k)}
\renewcommand\sp{\operatorname{Spec}}
\renewcommand\sf{\operatorname{Spf}}
\newcommand\proj{\operatorname{Proj}}
\newcommand\aut{\operatorname{Aut}}
\newcommand\grv{{\operatorname{Gr}}(V)}
\newcommand\grass{\operatorname{Grass}}
\newcommand\gr{\operatorname{Gr}}
\newcommand\grb{{\operatorname{Gr}}^\bullet}
\newcommand\glv{{\operatorname{Gl}}(V)}
\newcommand\glve{{\widetilde{\operatorname{Gl}}(V)}}
\newcommand\gl{\operatorname{Gl}}
\newcommand\pic{\operatorname{Pic}}
\renewcommand\hom{\operatorname{Hom}}
\renewcommand\det{\operatorname{Det}}
\newcommand\detd{\operatorname{Det}^\ast}
\newcommand\im{\operatorname{Im}}
\newcommand\res{\operatorname{Res}}
\newcommand\limi{\varinjlim}
\newcommand\limil[1]{\underset{#1}\varinjlim\,}
\newcommand\limp{\varprojlim}
\newcommand\limpl[1]{\underset{#1}\varprojlim\,}

\renewcommand\o{{\mathcal O}}
\renewcommand\L{{\mathcal L}}
\renewcommand\c{{{\mathcal C}}^\bullet}
\renewcommand\P{{\mathbb P}}
\newcommand\Z{{\mathbb Z}}
\newcommand\N{{\mathbb N}}
\newcommand\A{{\mathbb A}}
\newcommand\M{{\mathcal M}}
\newcommand\s{{\mathcal S}}
\newcommand\B{{\mathcal B}}

\newcommand\w{\widehat}
\renewcommand\tilde{\widetilde}

\newcommand\gen[1]{<\!#1\!>}
\newcommand\iso{@>{\sim}>>}
\renewcommand\lim{\underset{\underset{ A\sim V^+}A}\varprojlim}
\newcommand\fu{\underline}
\newcommand\beq{
      \setcounter{equation}{\value{defn}}\addtocounter{defn}1
      \begin{equation}}

\newcommand\dossub[2]{\underset{{\tiny\begin{gathered}
\raisebox{-4pt}{$#1$} \\ \raisebox{3pt}{$ #2$}
\end{gathered}}}}

\renewcommand{\thesubsection}{\thesection.\Alph{subsection}}

\title[Grassmannian of $k((z))$]{Grassmannian of $k((z))$:
Picard Group, Equations and Automorphisms.}

\author[F. J. Plaza Mart\'{\i}n]
{Francisco J. Plaza Mart\'{\i}n
\\
\smallskip\tiny Departamento de Matem\'aticas \  \  --\  \
Universidad de Salamanca
\\{\bf Postal Address:} Plaza de la Merced 1-4 \\
 37008 Salamanca -- SPAIN}

\thanks{
This work is partially supported by the
CICYT research contract n. PB96-1305.}

\email{fplaza@@gugu.usal.es}

\begin{abstract}
This paper aims at generalizing some geometric
properties of Grassmannians of finite dimensional vector spaces to
the case of Grassmannnians of infinite dimensional ones, in
particular for $\gr(k((z)))$. It is shown that the Determinant Line
Bundle generates its Picard Group and that the Pl\"ucker equations
define it as closed subscheme of a infinite projective space.
Finally,  a characterization of finite dimensional projective spaces
in Grassmannians allows us to offer an approach to the study of the
automorphism group.
\end{abstract}

\maketitle


\setcounter{tocdepth}1
\tableofcontents

\section{Introduction}

Infinite Grassmannians have recently emerged as a powerful
tool to study certain moduli spaces as well as being fundamental
objects in some algebraic approaches to Conformal Field Theories and
String Theory. It seems reasonable that a better understanding of
such spaces could throw some light on all these problems and the
pursuit of this motivates the present work.

Before explaining how the paper is organized, we wish to remark that
the method followed to reach the above-mentioned results is
essentially to relate  the infinite dimensional Grassmannian and a
suitable increasing chain of finite dimensional Grassmannians
contained in it (basically given in \cite{SW}).

Section \S2 is introductory and should help us to fix notations and
become familiarized with the main object of our study: infinite
Grassmannian. Certain known facts (\cite{AMP,BS,MP}) on infinite
Grassmannians are addressed such as their scheme structure and the
construction of global sections of (the dual of) the determinant
bundle, $\detd_V$. The rest of the section contains a deep study of
the relationships of finite and infinite Grassmannians, in particular
the restriction of sections of $\detd_V$, as well as some auxiliary
schemes such as the Grassmannian of the ``dual'' space and of a
metric space.

In section \S3 it is proved that the Picard group of (the index $0$
connected component of) the infinite Grassmannian is $\Z$ and that
the determinant bundle generates it. These results follow from
study of the restriction homomorphisms (from the infinite
Grassmannian to a finite one) induced on the group of Cartier
divisors.

The first part of section \S4 proves that the Pl\"ucker morphism
is a closed immersion. In the second
part, equations defining $\gr(k((z)))$ as a subscheme of $\P\Omega^*$
are computed.

The last section, \S5, begins with a characterization of finite
dimensional projective spaces in (finite or infinite)
Grassmannians (Theorem \ref{thm:proj-sp}). This result is the
cornerstone of Corollary \ref{cor:aut-fin-grass} which recovers all
the known results about the automorphism group of Grassmannians. In
order to deal with the infinite dimensional case, some extra
structure is added to $V$. It is now shown that the group of
automorphisms preserving this structure is a group deeply related to
the linear group.

I wish to express my gratitude to Prof. J.M. Mu\~noz Porras for his
supervision  and for the fruitful discussions we have had.

\section{Backgrounds on Grassmannians}

\subsection{Infinite Grassmannian}\label{subsec:inf-grass}

In \cite{PS,SW} the infinite Grassmannian of $V$, a Hilbert space
over a field $k$, and a certain decomposition of it, $V\simeq
V^-\oplus V^+$, is defined as the set of subspaces $L\subseteq V$
such that the projection $L\to V^-$ has finite dimensional kernel and
cokernel (our $+/-$ notation and that of Segal-Wilson have been
exchanged).

This definition of the above infinite Grassmannian, however, can be
slightly weakened (\cite{AMP,BS}). Consider the following set:
\beq\left\{\begin{gathered}\text{subspaces $L\subset V$ such that
$L\cap V^+$}\\
\text{and $V/(L+ V^+)$ are finite dimensional}
\end{gathered}\right\}
\label{grpoint}
\end{equation}
It then depends only on the choice of a subspace $V^+$ (instead of on
a decomposition). But let us give an algebraic definition of it.

We shall fix a pair $(V,\B)$ of a vector space and a family
of neigbourhoods of $(0)$ such that:
\begin{enumerate}
\item $A,B\in\B\quad \implies\quad A+B,A\cap B\in\B$,
\item $A,B\in\B\quad \implies\quad \dim(A+B)/A\cap B <\infty$ (that
is, $A$ is {\it commensurable} with $B$ and will be denoted by
$A\sim B$),
\item the topology is separated: $\cap_{A\in\B} A=(0)$,
\item $V$ is complete: $V\to\limpl{A\in\B}V/A$,
\item every finite dimensional subspace of $V/B$ is a
neigbourhood of $(0)$: $\limil{A\in\B}A\to V/B$ is surjective (for
$B\in\B$).
\end{enumerate}

The completion with respect to the topology induced by $\B$
will be denoted with $\;\hat{}\;$. If $T\to S$ is a morphism of
schemes over $k$ and
$U\subseteq V_S$ is a sub-$\o_S$-module, then $\hat U_T$ will denote
the $\o_T$-module $U\hat\otimes_{\o_S}\o_T$ (completion with respect
to the topology induced in $U$ by the restriction of the
natural topology of $V_S$).

The important fact, is that the set ({\ref{grpoint}}) is actually the
set of rational points of a reduced and separated scheme,
$\grb(V)$, over $\spk$ (as was shown in \cite{AMP} and in
$\cite{BS}$ \S4.3).  There it is proved that the
affine schemes
$F_A:=\hom_k(L_A,A)$ (where $A\in \B$ and $L_A\oplus A\simeq V$) are
an open covering of $\grb(V)$. Note that $F_A$ represents the functor:
$$S\rightsquigarrow
\left\{\begin{gathered} L\subseteq\hat V_S\text{ quasi-coherent} \\
\text{such that }L\oplus\hat A_S\simeq \hat V_S\end{gathered}\right\}
$$

This idea is essentially the same as in the case of Grassmannians of
finite dimensional vector spaces (see
\cite{EGA} I.9.7). But giving the functor of points of $\grb(V)$
is by no means trivial, since one should define the set of $S$-valued
points of $\grb(V)$ for an arbitrary scheme $S$ (\cite{AMP}).

From now on we will fix a subspace $V^+\in\B$. This choice allows to
introduce the index function as well as the determinant line bundle.
Let $\B_0$ be the set of subspaces $A\in\B$ such that $\dim A/(A\cap
V^+)=\dim V^+/(A\cap V^+)$.

Note that the index function:
$$\begin{aligned}\grb(V) &@>i>> \Z\\
L&\mapsto \dim_k(L\cap V^+)-\dim_k(V/L+V^+)
\end{aligned}$$
is locally constant, and moreover $\grb(V)=\coprod \gr^n(V)$,
where
$\gr^n(V):=i^{-1}(n)$ is connected. For simplicity's sake,
$\gr^0(V)$ will be denoted by $\grv$.

Over the scheme $\grv$ there is a canonical line bundle: the
determinant bundle. In order to define it, one considers $\L$, the
universal object of $\grv$. It is now not difficult to show that
the complex $\L\to \pi^*(V/V^+)$ ($\pi:\grv\to\spk$) of
$\o_{\grv}$-modules is a perfect complex (\cite{AMP}). We know from
\cite{KM} that its determinant, $\det_V$, exists.

Moreover, the complexes $\c_A\equiv\L @>{\pi_A}>>\pi^*(V/A)$ (for
$A\in\B_0$) are all canonically quasi-isomorphic, and
their determinants are therefore isomorphic in a canonical way (since
$\operatorname{End(\det_V)}=H^0(\grv,\o_{\grv})=k$, see \cite{AMP}).
Nevertheless, this allows us to construct global
sections of
$\detd_V$, since the morphism $\pi_A$ induces a section $det(\pi_A)\in
H^0(\grv,\detd(\c_A))$, and hence a section $\Omega_A\in
H^0(\grv,\detd_V)$, via the isomorphism
$\c_A\iso\c_{V^+}$. Note that $\Omega_A$ does not vanish at a point
$L\in\grv$ if and only if $\pi_A\vert_L:L\to V/A$ (the
restriction of $\c_A$ to the point $L$) is an isomorphism, i.e.
$L\oplus A\simeq V$.

Assume that it is possible to give  isomorphisms
$\phi_{AB}:\det(\c_A)\iso\det(\c_B)$ in a compatible way; that
is:
\beq
\phi_{AC}=\phi_{BC}\circ\phi_{AB}
\label{eq:compat}\end{equation}
We then have many sections:
$$\{\phi^*_{AV^+}(det(\pi_A))\,\vert\,A\in\B_0\}\,\subseteq\,
H^0(\grv,\detd_V)$$
that will enable us to introduce the Pl\"ucker morphism, since
$\{F_A\,\vert\,A\in\B_0\}$ is an open covering. If $\Omega\subseteq
H^0(\grv,\detd_V)$ is a subspace containing $\Omega_A$ for a
subcovering of $\{F_A\vert A\in\B_0 \}$, there is a canonical
surjection:
$$\Omega\otimes_k\o_\grv\to \detd_V$$
that by the universal property of the projective space induces a
morphism:
$$\grv\to\P\Omega^*$$
 which will be called the Pl\"ucker
morphism. (Here and in the sequel
$\P E^*$ will  denote the scheme $\proj(S^\bullet E)$, for a
$k$-vector space $E$).

When studying the Grassmannian of $k((z))$ we assume that:
$$\B\,:=\,\left\{\begin{gathered}
\text{subspaces $A\subseteq V$ containing $z^n\cdot k[[z]]$ as a}\\
\text{subspace of finite codimension (for $n\in{\mathbb Z}$)}
\end{gathered}\right\}$$
and $V^+:=k[[z]]$.

Let us now compute sections  explicitly. Let us denote with
$\s$ the set of Maya diagrams of virtual cardinal zero; that is, the
set of strictly increasing sequences of integers
$S=\{s_i\}_{i\geq 0}$ such that:
\begin{itemize}
\item there exists an integer $i_0$ such that
$\{i_0,i_0+1,\dots\}\subseteq S$,
\item $\#(S\cap\Z_{<0})=\#(\Z_{\geq 0}-S)$ (condition of
virtual cardinal zero).
\end{itemize}

Let $e_i$ be $z^i\in V$, and let
$A_S$ be $z$-adic completion of the subspace $\gen{\{e_{s_i}\}_{i\geq
0}}$ for a given Maya diagram $S$. Observe that for $A\in\B $ there
exists a Maya diagram such that $ A+ V^+\subseteq  A_S$. Moreover, if
$L\in F_A(k)\subset\grv$ then there exists an $S'\in\s$ (contained in
$S$) such that $L\in F_{A_{S'}}(k)$; that is, the open subschemes
$\{F_{A_S}\vert S\in\s\}$ cover $\grv$.

Fix a Maya diagram $S=\{s_i\}\in\s$. Following
\cite{AMP} one obtains a canonical isomorphism:
{\small \beq\detd\c_{A_S}\otimes \wedge(A_S/A_S\cap
V^+)\otimes\wedge(V^+/A_S\cap V^+)^*\,\iso\,\detd_V
\label{det-isom}\end{equation}}
(where $\wedge$ denotes the highest exterior power).

Since we have fixed a ``dense'' basis $\{e_i\}$ of $V$, we can give
a generator of $\wedge(A_S/A_S\cap V^+)\otimes\wedge(V^+/A_S\cap
V^+)^*$. Let $J=\{j_1,\dots,j_n\}=(\Z_{<0}\cap S)$ and
$K=\{k_1,\dots,k_n\}=(\Z_{\geq 0}- S)$. Then:
{\small \beq<e_J\otimes e_K^*>\,=\,\wedge(A_S/A_S\cap
V^+)\otimes\wedge(V^+/A_S\cap V^+)^*
\label{ejek}\end{equation}}
\noindent where $e_J:=e_{j_1}\wedge\dots\wedge e_{j_n}$ and
$e^*_K:=e^*_{k_1}\wedge\dots\wedge e^*_{k_n}$ (where
$e^*_i(e_j)=\delta_{ij}$). Moreover, from the
choice of a basis of $V$ in order to obtain the isomorphisms
({\ref{det-isom}}) it follows that the compatibility condition
({\ref{eq:compat}}) is satisfied. Summing up, if $s$ is a global
section of $\detd\c_{A_S}$ then (by ({\ref{det-isom}}) and
({\ref{ejek}})) $s\otimes e_J\otimes e_K^*$ is a global section of
$\detd_V$.

Let us denote with $\Omega_S$ the section associated
with $A_S$. In the sequel $\Omega$ will denote the vector space
$\gen{\{\Omega_S\}_{S\in\s}}$., and observe that the Pl\"ucker
morphism $\grv @>{\frak p}>> \P\Omega^*$ is well defined.

\subsection{Finite Grassmannians}\label{subsect:fin}

Let us now assume that $V$ is a finite dimensional $k$-vector space
and that $\B$ is the set of subspaces of $V$. Choose a basis
$\{e_1,\dots,e_d\}$ of
$V$ such that
$V^+=\gen{e_1,\dots,e_{d-r}}$, and let
$\{e^*_1,\dots,e^*_d\}$ be its dual basis. Let $\s$ now be the
set of strictly increasing sequences of $d-r$ integers $S\equiv
0<s_1<\cdots< s_{d-r}\leq d$; and for $S\in\s$ define
$A_S=\gen{e_{s_1},\dots,e_{s_{d-r}}}$.

Observe that $\{F_{A_S}\vert S\in\s\}$ is again a covering of
$\grv$ and that the rational points of $\grv$ are precisely
the $r$-dimensional subspaces.

What is $\Omega_S(L)$ now? ($L$ being a rational point of
$\grv$, and $S\in\s$). Note that by the basic properties of
$\det$ (\cite{KM}), one has:
$$\detd\c_{A_S}\vert_L\,\iso\, \wedge L^*\otimes\wedge (V/A_S)$$
and hence:
{\small\beq
\begin{aligned}\Omega_S(L) & =det(\pi_{A_S}\vert_L)\otimes e_J\otimes
e_K^*= \\
& = \pi_{A_S}(l_1)\wedge\dots\wedge\pi_{A_S}(l_r)\otimes
l_1^*\wedge\dots\wedge l_r^*\otimes e_J\otimes e_K^*
\end{aligned}
\label{omegal}\end{equation}}
where:
\begin{itemize}
\item $\{l_1,\dots,l_r\}$ is a basis of $L$ and
$\{l^*_1,\dots,l^*_r\}$ its dual (note that ({\ref{omegal}}) does not
depend on the choice of the basis),
\item $e_J=e_{j_1}\wedge\dots\wedge e_{j_n}$ where
$\{e_{j_1},\dots, e_{j_n}\}$ generates $A_S/A_S\cap V^+$; that is,
$\{j_1,\dots,j_n\}=\{s_1,\dots,s_{d-r}\}-\{1,2,\dots,d-r\}$,
\item $e_K^*=e^*_{k_1}\wedge\dots\wedge e^*_{k_n}$ where
$\{e^*_{k_1},\dots, e^*_{k_n}\}$ generates $(V^+/A_S\cap V^+)^*$; that
is,
$\{k_1,\dots,k_n\}=\{1,2,\dots,d-r\}-\{s_1,\dots,s_{d-r}\}$.
\end{itemize}
Computing ({\ref{omegal}}) one has that:
{\small $$\Omega_S(L)\,=\,
(e^*_{\bar s_1}\wedge\dots\wedge e^*_{\bar  s_r})
(l_1\wedge\dots\wedge l_r)
\cdot e_1\wedge\dots\wedge e_{d-r}\otimes
l_1^*\wedge\dots\wedge l_r^*\,\in\,(\detd_V)\vert_L$$}
where $\{\bar s_1,\dots,\bar
s_r\}:=\{1,\dots,d\}-\{s_1,\dots,s_{d-r}\}$; that is, $\{e_{\bar
s_1}^*,\dots,e_{\bar s_r}^*\}$ generates $(V/A_S)^*$ and
$\{l_1,\dots,l_r\}$ is a basis of $L$.

Using this calculations, one deduces a natural isomorphism:
$$\begin{aligned}
H^0(\grv,\detd_V)&\to \wedge^rV^* \\
\Omega_S &\mapsto e^*_{\bar s_1}\wedge\dots\wedge e^*_{\bar  s_r}
\end{aligned}$$
where $\{\bar s_1,\dots,\bar s_r\}=\{1,\dots,d\}-
\{s_1,\dots,s_{d-r}\}$.

\subsection{Morphisms between Grassmannians}\label{subsec:morph}

Let us now relate both Grassmannians, the infinite one and the
finite one. Let $(V,\B,V^+)$ be as usual, and let $N$ be a rational
point of $\grb(V)$ and $M\subseteq V$ a subspace such that
$N\subset M$. If $\grass(M/N)$ is the standard Grassmannian of a
vector space (see \cite{EGA}~{\bf I}) we know from \cite{AMP} that the
morphism:
$$\begin{aligned}
j:\grass(M/N) &\to \grb(V) \\
L&\mapsto \pi^{-1}(L) \end{aligned}$$
(where $\pi:M\to M/N$) is a closed immersion. Moreover, if
$M\in\grb(V)$, there exists a canonical isomorphism:
$$j^*\det_V\iso \wedge\L^f\otimes\wedge(M/N)^*\otimes
\wedge(V^+\cap M)\otimes\wedge(V/M+V^+)^*$$
where $\L^f$ is the universal submodule of $\grass(M/N)$.
A very important consequence of this fact, is that the
morphism:
$$\grass^k(M/N) @>j>> \grv @>{\frak p}>> \P H^0(\grv,\detd_V)^*$$
factors through the Pl\"ucker morphism of $\grass^k(M/N)$ (for a
certain $k\in\Z$).

Looking at restrictions, we note first that if $A\in\B $ satisfies
$V/M+A=(0)$ and $N\cap A=(0)$, then $j^{-1}(F_A)=F_{M\cap A}$. It
follows that the restriction of the section $\Omega_A$ is
$\Omega_{M\cap A}$. One now concludes that the restriction
homomorphism of global sections:
\beq
\Omega \to H^0(\grass^k(M/N),\detd_{M/N})
\label{eq:restr-surj}\end{equation}
is surjective.

In the special case of $\gr(k((z)))$, there exist
Maya diagrams $S(M)$ and $S(N)$ (but not of virtual cardinal zero)
such that $M\in F_{A_{S(M)}}$ and $N\in F_{A_{S(N)}}$, where we can
assume that  $S(M)\subseteq S(N)$, since $N\subseteq M$. The
restriction of $\Omega_S$ for $S\in\s$ is now:
$$j^*\Omega_S\,=\,\cases \Omega_{\tilde S} \text{ where $\tilde
S:=S-S(N)$}& \text{ if }S(M)\subseteq S\subseteq S(N) \\
0 &\text{otherwise}\endcases$$

The case of a morphism between two finite Grassmannians is quite
similar. Let $\pi:V\to V'$ be a surjective
morphism between two $k$-vector spaces of
dimensions $d$ and $d'$ respectively, and let
$\bar d$ be $\dim_k(\ker\pi)=d-d'$. There is a natural
morphism:
$$\begin{aligned}
\grass^r(V') &\to \grass^{r+\bar d}(V) \\
L&\mapsto \pi^{-1}(L) \end{aligned}$$
which is known to preserve the determinant sheaf, and it
therefore induces a restriction homomorphism between global sections:
\beq\wedge^{r+\bar d}V^*\to\wedge^r{V'}^*
\label{fin-rest}\end{equation}
given by the inner contraction with $\wedge\ker\pi$. Choose a basis
$\{e_1,\dots,e_d\}$ of $V$ and $\{e'_{\bar d+1},\dots,e'_d\}$ of
$V'$, such that:
$\pi(e_i)=0$ for $i\leq\bar d$, and $\pi(e_i)=e'_i$ for $i>\bar d$.
The homomorphism ({\ref{fin-rest}}) is now written as:
$$e^*_{j_1}\wedge\dots\wedge e^*_{j_{r+\bar d}}\,\mapsto\,
\cases {e'}^*_{j_{1+\bar d}}\wedge\dots\wedge {e'}^*_{j_{r+\bar
d}} & \text{ if $j_i=i$ for $i\leq \bar d$,} \\
0 &\text{ otherwise.}\endcases$$
(where we assume that $1\leq j_1<j_2<\dots<j_{r+\bar d}\leq d$).

\subsection{Related Grassmannians}

\subsubsection{The Grassmannian of the dual space}
Let $(V,\B,V^+)$ be as usual. For a
given submodule $U\subseteq \hat V_S$ ($S$ a $k$-scheme), we
introduce the following notation:
$$\begin{aligned}
U^* &\,:=\, \hom_{\o_S}(U,\o_S) \\
U^c &\,:=\, \{f\in U^* \text{ continuous }\}
\end{aligned}$$
 where the topology in $U$ is given by $\{\hat A_S\cap U\vert A\in\B
\}$, and $\o_S$ has the discrete topology. And define:
$$\begin{aligned}
U^\circ &\,:=\, \{f\in (\hat V_S)^* \,\vert\, f\vert_U\equiv 0 \} \\
U^\diamond &\,:=\, \{f\in (\hat V_S)^c \,\vert\, f\vert_U\equiv 0 \}
\end{aligned}$$

In order to make explicit the meaning of the expression
``Grassmannian of the dual (continuous) space'' we consider in
$V^c$ the family:
$$\B^\diamond \,:=\,\{A^\diamond\text{ where }A\in\B\}$$

An easy consequence of linear algebra is the following:

\begin{lem}\hfill
\begin{enumerate}
\item $V^c=\limi A^\diamond$ ($A\in\B$),
\item $V=\limi A$ ($A\in\B$),
\item the topology induced by $\B^\diamond$ is separated,
\item $V^c$ is complete.
\end{enumerate}
\end{lem}

It is now easy to check that the family $\B^\diamond$ satisfies the
same conditions as $\B$ (given at the begining of
\ref{subsec:inf-grass}), and hence one obtains the following:

\begin{thm}
The Grassmannian of $(V^c,\B^\diamond)$ exists if and only if that of
$(V,\B)$ does.
\end{thm}

Here, we shall construct a canonical
isomorphism between the Grassmannian of $(V,\B)$ and that of
$(V^c,\B^\diamond)$. The expression of this isomorphism for the
rational points will be that given by the incidence:
\beq
\begin{aligned} I:\grb(V) &\longrightarrow \grb(V^c)
\\ L\quad &\longmapsto\quad L^\diamond \end{aligned}
\label{eq:inci}\end{equation}

However, the existence of such a morphism is equivalent to show that
$\L^\diamond$ is in fact a point of $\grb(V^c)$. Note that
$\L^\diamond$ is a quasi-coherent sheaf since it is a direct limit of
quasi-coherent sheaves; namely:
$$\L^\diamond\,=\,\limil{A\in\B } \big( V/A+\L)^*$$
Since $\{F_A\}_{A\in\B }$  is a covering of
$\grb(V)$, it is thus enough to prove that $\L^\diamond$ and $\hat
A^\diamond_{F_A}$ are in direct sum (over $F_A$), but this follows
from the canonical isomorphism $\hat V_{F_A}\simeq
\L_{F_A}\oplus\hat A_{F_A}$.

An easy calculation shows now the relationships between the
determinant bundles and index functions of both Grassmannians. The
following formulae hold:
$$I^*\det_{V^c}\simeq \detd_V\qquad\qquad I^*(i_{V^c})=-i_V$$
 ($i$ being the index function).

\subsubsection{The case of a metric space}\label{subsubsect:metric}

Assume now that there is a given irreducible and hemisymmetric metric
on $V$:
$$T_2:V\times V\to k$$

A subspace $L\subseteq V$ is called totally isotropic when
$L\subseteq L^\perp$. Observe that since $T_2$ is hemisymmetric the
condition of maximal totally isotropic (m.t.i.) on a subspace $L$
implies $L=L^\perp$. The subspace $V^+$ will be assumed to be m.t.i.

It is not difficult to prove that under these hypothesis the polarity:
$$iT_2:V\iso V^c$$
is a bicontinuous isomorphism of vector spaces (with respect to the
topologies induced by $\B$ in $V$ and by $\B^\diamond$ in $V^c$) and
therefore it induces an isomorphism:
$$\grb(V^c)\,\iso\,\grb(V)$$

The composition of the latter isomorphism and ({\ref{eq:inci}}) is  a
natural involution of the Grassmannian, whose expression at the
rational points is given by:
\beq
\aligned R:\grb(V) &\longrightarrow\grb(V) \\
L &\longmapsto L^{\perp} \endaligned
\label{eq:Raut}\end{equation}

Trivial calculation shows that $R^*\det_V\simeq\det_V$ and that
the index of a point $L\in\grb(V)$ is exactly the opposite of the
index of $R(L)=L^\perp$.

\section{Picard Group of $\grv$}

\begin{thm}
Let $k$ be an algebraically closed field. Assume that
$\dim(\grv)\geq 1$.
Then, the Picard group of $\grv$ is isomorphic to $\Z$ and the line
bundle $\det_V$ is a generator.
\end{thm}

\begin{pf}
Let us first deal with the $\dim(\grv)=1$ case. Recall that
$F_A\iso\hom(L, A)$ (where $V\simeq L\oplus A$ and $A\in\B$)  is
an open subscheme of $\grv$ and that its dimension equals
$\dim_kL\cdot\dim_k  A$. It follows that $\dim_kL=\dim_k A=1$
and that $\grv$ is the projectivization of a $2$-dimensional vector
space; i.e. $\grv$ is the projective line, whose Picard group is
$\Z$, as is well known.

For the general case, recall that the Picard group of $\grv$ is
canonically isomorphic to the class group of Cartier divisors, since
$\grv$ is integral.

Fix $A\in\B_0$ and assume $Z_A:=\grv-F_A$ (the locus where $\Omega_A$
vanishes) to be irreducible. It then implies the exactness of the
sequence:
$$\begin{aligned}
\Z &\to \pic(\grv)\to \pic(F_A)
\\ 1 &\mapsto \o_{\grv}(-Z_A)=\det_V
\end{aligned}$$
Observe that $\pic(F_A)=0$, since $F_A$ is the spectrum of
a ring of polynomials (in finitely or infinitely many variables)
which is factorial (see \cite{Bou} cap. VII \S3.5.). Finally, the
line bundle ${\detd_V}^{\otimes n}$ cannot be trivial for $n>0$,
because its space of global sections has dimension greater than $1$
and $H^0(\grv,\o_{\grv})$ has dimension $1$ (\cite{AMP}).

The proof is therefore reduced to the following:
\end{pf}

\begin{lem}
Let $k$ be an algebraically closed field, and $\dim(\grv)\geq 2$.
There exists a subspace $A\in\B_0$ of $V$ such that
$Z_A$ is irreducible.
\end{lem}

\begin{pf}
If $\grv$ is of finite type (that is, $V$ is finite dimensional)
 the proof is an easy consequence of the Bertini Theorem (see
\cite{Ha} II.8.18 and III.7.9.1), since $Z_A$ is precisely a
hyperplane section of the Pl\"ucker morphism. Moreover, it implies
that $Z_A$ is irreducible for $A$ generic.

Assume now that $V$ is not finite dimensional. We shall prove that if
$Z_A$ is reducible then its natural restriction to a (suitable)
finite dimensional Grassmannian is also reducible. Bearing this in
mind and since the restriction homomorphism of global
sections is surjective ({\ref{eq:restr-surj}}), one concludes that
this is not possible for generic $A$, and the result follows.

Let $Z_1,Z_2\subset Z_A$ be two (different) irreducible components
of $Z_A$. If $Z_1\cap Z_2\neq\emptyset$, take a point $L_0\in Z_1\cap
Z_2$ and a subspace $B\in\B_0$ such that $L_0\in F_B$. Let $L_1,L_2$
be points of $Z_1\cap F_B$, $Z_2\cap F_B$ respectively, such that
$\dim_k(L_i/L_0\cap L_1\cap L_2)<\infty$ (for $i=1,2$). Note that
these three points can be assumed to be rational. From Subsection
{\ref{subsec:morph}} we know that a certain connected component,
$G$, of $\grass((L_0+L_1+L_2)/(L_0\cap L_1\cap
L_2))$ is naturally mapped into $\grv$.
Obviously, the points $L_0,L_1$ and  $L_2$ lie in $G$, and
$\emptyset\neq G\cap Z_i\neq G$ ($i=1,2$) and $G\cap Z_1\neq G\cap
Z_2$. It follows that the divisor $G\cap Z_A\subset G$ is a reducible
hyperplane section of $G$.

Let us finally prove the other case: $Z_A$ is reducible and
$Z_1,Z_2\subset Z_A$ are two irreducible components such that
$Z_1\cap Z_2=\emptyset$. Consider two rational points  $L_i\in Z_i$
($i=1,2$) such that $\dim_k(L_i/L_1\cap L_2)<\infty$, and let $G$
be the connected component of  $\grass((L_1+L_2)/(L_1\cap L_2))$ whose
image lies in $\grv$. One then concludes using similar ideas as
before.
\end{pf}

\section{Pl\"ucker equations}

\subsubsection*{The Pl\"ucker morphism}

We shall show that the Pl\"ucker morphism is a closed immersion.
The proof is based on that of \cite{EGA}~{\bf I}.9.8.4 for the
finite dimensional case.

From \cite{AMP} one canonically obtains a
$1$-dimensional subspace, $<\Omega_A>$ ($A\in\B_0$), of
$H^0(\grv,\detd_V)$ by using the canonical isomorphism:
{\small\beq
\det^*\c_A\otimes (\wedge A/A\cap V^+)\otimes (\wedge V^+/A\cap
V^+)^*\iso \detd_V
\label{eq:det-iso}\end{equation}}
 and the canonical section of $\det^*\c_A$. Let
$\Omega$ be the subspace:
$$\sum_{A\in\B_0 }<\Omega_A>\subseteq H^0(\grv,\detd_V)$$

But let us offer another description of this. Observe that $\Omega$
is the image of the morphism of vector spaces:
$$\Psi:\bigoplus_{A\in\B_0} (\wedge A/A\cap V^+)\otimes (\wedge
V^+/A\cap V^+)^* \,\to\, H^0(\grv,\detd_V)$$
defined by tensoring the $A$-component with $det(\pi_A)$.

\begin{lem}\label{lem:omega}
There exists a natural surjective morphism:
$$\CD
\underset{A\in\B_0}\bigoplus (\wedge A/A\cap V^+)\otimes
(\wedge V^+/A\cap V^+)^*
\\  @V{\Phi}VV \\
\dossub{B\in\B}{B\subseteq V^+}\limi \left(
\bigoplus_{k=0}^{\dim(V^+/B)}  \wedge^k V/{V^+}
\otimes
\big(\wedge^k V^+/B\big)^*\right)
\endCD$$
such that $\Psi$ factors through $\Phi$ and an injection:
{\small $$\dossub{B\in\B}{B\subseteq V^+}\limi
\left(\bigoplus_{k=0}^{\dim(V^+/B)} \wedge^k V/{ V^+}
\otimes\big(\wedge^k V^+/B\big)^*\right)
\hookrightarrow H^0(\grv,\detd_V)$$}
\end{lem}

\begin{pf}
Let us fix $B\in\B $ such that $B\subseteq V^+$. Note that the
morphism:
$$\dossub{A\in\B_0}{A\cap V^+=B}\bigoplus (\wedge
A/A\cap V^+)\otimes (\wedge V^+/A\cap V^+)^* \,\to\,
H^0(\grv,\detd_V)$$
factors through  a surjection onto $\wedge^k
V/{ V^+} \otimes (\wedge^k V^+/B)^*$, where $k=\dim_k(V^+/B)$.
Therefore, the linear map
$\Psi$ factors through a surjection onto:
$$\bigoplus_{B\subseteq V^+} \left(\wedge^k V/{ V^+}\otimes
\big(\wedge^k V^+/B\big)^*\right)$$

Fixing $B$ again, observe further that the induced morphism:
$$\dossub{B\subseteq B'}{B'\subseteq V^+}\bigoplus
\wedge^k V/{V^+} \otimes (\wedge^k V^+/B')^*
\,\to\, H^0(\grv,\detd_V)$$
($k=\dim_k(V^+/B)$) factors through a surjection onto:
\beq\bigoplus_{k=0}^{\dim(V^+/B)}  \wedge^k V/{ V^+}
\otimes\big(\wedge^k V^+/B\big)^*
\label{eq:dir-syst}\end{equation}

Note that the set ({\ref{eq:dir-syst}}) is a direct system as $B$
varies. The statement now follows easily.
\end{pf}

\begin{thm}\label{thm:plucker-imm}
The Pl\"ucker morphism:
$$\grv\longrightarrow \P\Omega^*$$
is a closed immersion.
\end{thm}

\begin{pf}
Let $U_A$ the affine open
subscheme of $\P\Omega^*$ where the $A$-coordinate has no zeroes.
It is clear that ${\frak p}(F_A)\subseteq U_A$, and that it is
 enough to see that ${\frak p}\vert_{F_A}$ is a closed
immersion.

For the sake of clarity, it will be assumed that $A=V^+$.
Nevertheless, the general case presents no extra difficulty.
By fixing sections of $V\to V/V^+$ and $\Omega\to
\Omega/\gen{\Omega_+}$ one has  identifications $F_{V^+}\simeq
\hom(V/ V^+, V^+)$ and $U_{V^+}\simeq
\hom(\Omega/\gen{\Omega_+},\gen{\Omega_+})$. The restriction of the
Pl\"ucker morphism to $F_{V^+}$ is now a morphism:
\beq\hom(V/  V^+, V^+) \to
\hom(\Omega/\gen{\Omega_+},\gen{\Omega_+})
\label{eq:plucker-res}\end{equation}

By Lemma {\ref{lem:omega}} one has that:
$$\Omega\,\simeq\, \dossub{B\subseteq V^+}{B\in\B}\limi
\left(\bigoplus_{k=0}^{\dim(V^+/B)} \wedge^k V/{V^+}
\otimes \big(\wedge^k V^+/B\big)^*\right)$$
and that $\gen{\Omega_+}$ corresponds to $B=V^+$. Recalling that
$ V^+=\limp V^+/B$, it
follows that ({\ref{eq:plucker-res}}) is the inverse limit of:
{\small $$\hom(V/ V^+,V^+/B)\to
\hom\left(\bigoplus_{k=1}^{\dim(V^+/B)} \wedge^k V/{ V^+}
\otimes\big(\wedge^k V^+/B\big)^*,\gen{\Omega_+}\right)$$}
where $B\in\B $ is such that $B\subsetneqq V^+$.

Observe that all
these spaces are affine schemes, and it then suffices prove that
given $B\in\B $ such that $B\subsetneqq V^+$ the morphism:
{\small $$\hom(V/ V^+,V^+/B)\to \prod_{k=1}^{\dim(V^+/B)}
\hom(\wedge^k V/{ V^+}
\otimes\big(\wedge^k V^+/B\big)^*,\gen{\Omega_+})$$}
is a closed immersion. This, however, is trivial since it is the
graph of a morphism.
\end{pf}

\subsubsection*{Restriction of sections}(Case of $V=k((z))$).

Take $L_i:=\gen{\{e_j\}_{j< i}}$ for an integer $i$. Note that:
$L_i\in\gr^\bullet(V) $. Denote by $j_i$ the induced morphism
$\grass(L_i/L_{-i})\to\gr^\bullet(V) $ and let $G_i$ be
$\grass(L_i/L_{-i})\cap\grv$. Observe that the index corresponding
to the connected component $G_i$ is $i$. The $\gr_0$ of the
Segal-Wilson paper (\cite{SW}) is precisely the union of all $G_i$.

Recall that the following diagram is commutative:
\beq\CD
G_i @>{\frak p}_i>> \P (\bigwedge^{i} L_i/L_{-i}) \\
@V{j_i}VV @VV{\iota_i}V \\
\grv @>{\frak p}>> \P \Omega^*
\endCD
\label{eq:plucker}\end{equation}
Denote by $\iota^*_i$ the restriction morphism $\Omega\to\
\bigwedge^{i}( L_i/L_{-i})^*$ which is known to be surjective
by ({\ref{eq:restr-surj}}). Let us denote by $\gen{X}$ the free
$k$-module generated by a set $X$. Let $\s_i$ be the set of strictly
increasing sequences $-i\leq s_0<s_1<\dots< s_{i-1}\leq i-1$. Note
that
$\bigwedge^{i}( L_i/L_{-i})^*\simeq <\s_i>$ (see
{\ref{subsect:fin}}),
$\Omega\simeq \gen{\s}$, and that $\iota^*_i$ is the morphism
induced by the map:
$$\begin{aligned}
\s & \longrightarrow \s_i \\
\{s_i\}_{i\geq 0} &\mapsto
\cases \{s_0,\dots,s_{i-1}\} & \text{if $-i\leq s_1$ and $s_j=j$ for
all
$j\geq i$}\\ 0 & \text{otherwise} \endcases
\end{aligned}$$
This morphism has a natural section; namely:
$$\begin{aligned}
\s_i & \longrightarrow \s \\
\{s_0,\dots,s_{i-1}\} &\mapsto \{s_0,\dots,s_{i-1},i,i+1,\dots\}
\end{aligned}$$
Let $\sigma_i$ be the induced morphism
$S^\bullet\gen{\s_i}\hookrightarrow S^\bullet\gen{\s}$ (where
$S^\bullet E$ denotes the symmetric algebra generated by a vector
space $E$).

In a similar way, and bearing in mind the restriction morphism of
global sections of determinant bundles on finite Grassmannians, one
constructs surjections $\s_j\to\s_i$ (for
$j\geq i$) and sections $\s_i\to\s_j$. It follows easily that the
family $\{\s_i\}_{i\geq 1}$ is a inverse system (with respect to the
surjections) and a direct system (with respect to the sections) in a
compatible way. Moreover, one has:
$$\s\simeq\limil{i}\s_i\,\iso\,\limpl{i}\s_i$$
from which one has that $S^\bullet\gen{\s}\iso \limi
S^\bullet\gen{\s_i}$; and hence:
\beq
I=\limi\big(I\cap \sigma_i(S^\bullet\gen{\s_i})\big)\iso \limi
\big(\iota_i^*I\cap S^\bullet\gen{\s_i}\big)
\label{eq:limi}\end{equation}
for every submodule $I\subset S^\bullet\gen{\s}$.

\subsubsection*{Pl\"ucker equations}(Case of $V=k((z))$).

We shall now give explicit equations for the infinite
Grassmannian, $\gr(k((z)))$, in a
infinite dimensional projective space. Such equations will in
fact be the infinite set of Pl\"ucker relations, and in this sense we
prove that our definition of infinite Grassmannian is the same as
that of Sato-Sato [{\bf SS}]. Further, in this setting the
Segal-Wilson Grassmannian $\gr_0$ (see \S2 of \cite{SW}) is to be
interpreted as the set of points of the Sato Grassmannian with
finitely many non-zero coordinates. In a certain sense, the scheme
$\gr(k((z)))$ unifies both Grassmannians.

By the compatibility of the Pl\"ucker morphisms of
({\ref{eq:plucker}}), one has also the following commutative diagram
of sheaves:
$$\CD
\big(S^\bullet \Omega \big){\widetilde{}}
@>{{\frak p}^*}>> \o_{\grv}
\\ @V{\iota^*_i}VV @VV{j_i^*}V \\
\big(S^\bullet (\bigwedge^{i} L_i/L_{-i})^*\big){\widetilde{}}
@>{{\frak p}_i^*}>>
\big(\underset{d\geq0}\oplus H^0(G_i,{\detd}^{\otimes
d})\big){\widetilde{}}
\endCD$$
where $\;\widetilde{}\;$ denotes the homogeneous localization sheaf.
Let
$I$ be the kernel of ${\frak p}^*$ and $I_i$ that of ${\frak
p}_i^*$. Theorem {\ref{thm:plucker-imm}} implies that these ideals are
the equations defining $\gr(k((z)))$ and
$G_i$, respectively.

From ({\ref{eq:limi}}), and since $S^\bullet (\bigwedge^{i}
L_i/L_{-i})^* \simeq S^\bullet\gen{\s_i}$, one has that
$I=\limi(\iota_i^* I_i)$. The same
argument also implies that $I$ is generated by its degree 2
homogeneous component, since $I_i$ does (this is a classical result).
That is, the ideal
$I$ is generated by the union of the generators of $I_i$ for all $i$.

Recall that a family of generators of $I_i$ is given by:
$$X_S\cdot X_{S'}-\sum_{l\geq0}X_{S_l}\cdot X_{S'_l} \qquad
\text{for }S,S'\in\s_i, k\geq 0$$
where:
$$\begin{aligned}
S_l\,:=\,& \{s_0,\dots,s_{k-1},s'_l,s_{k+1},\dots,s_{i-1}\}\\
S'_l\,:=\,&
\{s'_0,\dots,s'_{l-1},s_k,s'_{l+1},\dots,s'_{i-1}\}
\end{aligned}$$
for sequences:
$$\begin{aligned}
S\,\equiv\,& -i\leq s_0<s_1<\dots< s_{i-1}\leq i-1 \\
S'\,\equiv\,&  -i\leq s'_0<s'_1<\dots< s'_{i-1}\leq i-1
\end{aligned}$$
and $X_S$ is defined as follows:
\begin{enumerate}
\item for a sequence $S\in\s_i$, $X_S$ denotes the degree 1 element
of $S^\bullet\gen{\s_i}$ canonically  associated to $S$;
\item for an arbitrary set of $i$ distinct numbers
$\{s_0,\dots,s_{i-1}\}$ (such that $-i\leq s_j\leq i-1$ for all
$j$),  let $\sigma$ be the permutation of $i$-letters such that
it takes $S$ in increasing order, and then
$X_S:=\operatorname{sig}(\sigma)\cdot X_{S^\sigma}$;
\item $X_S:=0$ otherwise.
\end{enumerate}

We have thus proved the following:

\begin{thm}\label{thm:plucker}
The ideal $I$  is generated by its
degree $2$ homogeneous component. More precisely, it is
generated by the following equations:
$$X_S\cdot X_{S'}-\sum_{l\geq0}X_{S_l}\cdot X_{S'_l} \qquad
\text{for }S,S'\in\s, k\geq 0$$
where $X_{S_l}$ is defined as above.
\end{thm}

Note that these sums are
finite. These equations will be called the set of all Pl\"ucker
relations.

\begin{cor}
The closed subscheme of $\P\Omega^*$ defined by the Pl\"ucker
relations coincide with $\gr(k((z)))$ (via the Pl\"ucker morphism).
\end{cor}

\section{Automorphisms of Grassmannians}

Our goal now is to describe the relation between the linear
group of a vector space and the automorphism group of
its Grassmannian (in the infinite dimensional case). It is a classical
result that both groups coincide in the finite dimensional case.

\subsection{The Linear Group}\label{subsec:linear-group}

Fix a pair $(V,\B)$ as usual.
Given a $k$-scheme $S$, $\aut_{\o_S}(\hat V_S)$ will denote the
automorphism group of $\hat V_S$ as an $\o_S$-module.

\begin{defn}\hfill
\begin{itemize}
\item A sub-$\o_S$-module $A\subseteq \hat V_S$
belongs to $\B$ if there exists $B\in\B$ such that $\hat B_S\subseteq
A$ and the quotient is free of finite type.
\item An automorphism $g\in\aut_{\o_S}(\hat V_S)$ is bicontinuous
(w.r.t. $\B$) if there exists $A\in\B$ such that both $g(\hat A_S)$
and $g^{-1}(\hat A_S)$ belong to $\B$.
\item The linear group, $\glv$, associated to $(V,\B)$, is the
contravariant functor over the category of $k$-schemes given by:
$$S\rightsquigarrow \glv(S)=\{g\in\aut_{\o_S} (\hat
V_S)\text{ such that $g$ is bicontinuous }\}$$
\end{itemize}
\end{defn}

\begin{thm}\label{thm:detinv}
There exists a canonical action of $\glv$ on (the functor of points
of) $\grb(V)$:
$$\begin{aligned}
\glv\times & \grb(V) @>{\mu}>>\grb(V)  \\
(g, &L)\phantom{xx}\longmapsto g(L)
\end{aligned}$$
Moreover, this action preserves $\detd_V$.
\end{thm}

\begin{pf}
Fix $g\in\glv(S)$. Note that it suffices
check that $g(L)$ is a point of the Grassmannian for arbitrary
$L\in F_A(S)$.

From $L\oplus \hat A_S\simeq \hat V_S$ it follows that $g(L)\oplus
g(\hat A_S)\simeq \hat V_S$. Let $B\in\B$ be such that $g(\hat
A_S)/\hat B_S$ is free of finite type. It follows from \cite{AMP}
that $g(L)\cap\hat B_S=0$ and $\hat V_S/(g(L)+\hat B_S)$ is locally
free of finite type, and hence
$g(L)\in \grb(V)(S)$, as desired.

Let us explain what ``preserve'' means. Let $f:S\to \grv$ be an
$S$-valued point of $\grv$, $g$ an element of $\glv(S)$, and
$f_g:S\to \grv$ the transform of $f$ under $g$. The claim is:
$$f^*\det_V\,\simeq\, f_g^*\det_V$$
(as line bundles over $S$).
But this is a direct consequence of the properties of the determinant
(\cite{KM}) and the exactness of the sequence of complexes (written
vertically):
$$\CD
g(L)\oplus \hat B_S @>>> g(L)\oplus g(\hat V^+_S) @>>> g(\hat
V^+_S)/ \hat B_S  \\ @VVV @VVV @VVV \\
\hat V_S @>>>  \hat V_S @>>> 0
\endCD$$
(where $B\in\B$ is such that $g(\hat V^+_S)/\hat B_S)$ is free of
finite type).
\end{pf}

\subsection{Projective Spaces in $\grv$}\label{subsec:proj-sp}

For the sake of notation, let us denote simply by $D$ the dual of the
determinant sheaf $\detd_V$ over $\grv$, and by $D_L$ its stalk at a
rational point $L\in\grv$. Further, $\grv$ will be thought of as
a closed subscheme of $\P\Omega^*$ (Theorem {\ref{thm:plucker-imm}}),
and we shall make no distinction between a point of $\grv$ and the
same point considered as a point of $\P\Omega^*$. Let us first study
the structure of the projective lines contained in
$\grv$.

For a point $L\in\grv$, ${\cal H}_L$ is defined by:
$${\cal H}_L\,:=\,\left\{ \text{ hyperplanes
$\P\Omega^*$ passing through $L$ }\right\}$$
For a family $\{L_i\}$ of points of $\grv$ (from now on all the
points are assumed to be rational) we define:
$${\cal H}_{\{L_i\}}\,:=\, \bigcap_{i\in I}{\cal H}_{L_i}$$
Observe that (set-theoretically):
$${\cal H}_{\{L_i\}}\,=\,\text{Projectivization
of }\ker\left(\Omega\to
\bigoplus_i D_{L_i}\right)$$
Fix a family $\{L_i\}_{1\leq i\leq n}$, and assume  that there
exists $A\in\B $ such that $(V/( A+L_i))=(0)$. It then follows
that $\dim(A\cap L_i)$ is constant; call it $k$. Now, for $1\leq
j\leq n$ define:
$$\Lambda_j\,:=\, \bigwedge^k(A\cap L_1+\dots+A\cap L_j)^*$$

\begin{lem}
The restriction homomorphism:
$$\Omega\to \bigoplus_{i=1}^j D_{L_i}$$
factors  through $\Lambda_j$. Further,
$\Omega\to \Lambda_j$ is surjective.
\end{lem}

\begin{pf}
Recall that $D$ is isomorphic to the determinant of the
complex $\L\to V/A$ for $A\in\B$. Bearing in mind that there exists
an isomorphism $\wedge^k(A\cap L_i)^* \iso D_{L_i}$ (for $A$ as
above) and Lemma {\ref{lem:omega}}, the
claim reduces to an exercise of linear algebra.
\end{pf}

\begin{thm}
Three rational points $L_1,L_2,L_3$ of $\grv$ lie in a line (as
points of $\P\Omega^*$) if and only if $L_1\cap L_2\subseteq
L_3\subseteq L_1+L_2$ and both inclusions have codimension $1$.
\end{thm}

\begin{pf}
Observe that $\{L_1,L_2,L_3\}$ lie in a line if and only if:
\beq
{\cal H}_{L_1L_2L_3}\,=\,{\cal H}_{L_1L_2}
\label{eq:3line}\end{equation}

Consider the following
commutative diagram:
$$\CD
\Omega @>p_3>>\Lambda_3@>\rho_3>>\bigoplus_{i=1}^3 D_{L_i}
\\ @| @V{\pi}VV @V{\pi_{12}}VV \\
\Omega @>p_2>>\Lambda_2@>\rho_2>>\bigoplus_{i=1}^2 D_{L_i}
\endCD$$
and note that $p_3,p_2$ are surjective (by the above
Lemma), $\rho_2$ is also surjective (since $L_1$ and $L_2$ are
distinct), and finally $\pi$ is surjective too. Now, condition
({\ref{eq:3line}}) is equivalent to:
$$\ker(\rho_3)\,=\,\ker(\rho_2\circ\pi)$$
But from the very definition of $\Lambda_3$ one easily sees that:
$$\begin{gathered}
\ker(\rho_3)\,=\, (\wedge^k(A\cap L_1)+ \wedge^k(A\cap L_2)+
\wedge^k(A\cap L_3))^\circ
\,\subset\,\Lambda_3\\
\ker(\rho_2\circ\pi)\,=\, (\wedge^k(A\cap L_1)+ \wedge^k(A\cap
L_2))^\circ
\,\subset\,\Lambda_3\end{gathered}$$
and they are equal if and only if:
\beq
\wedge^k(A\cap L_3)\,\subseteq \, \wedge^k(A\cap L_1)+
\wedge^k(A\cap L_2)
\label{eq:inclu}\end{equation}

It is easy to see that ({\ref{eq:inclu}}) implies
that $A\cap L_3\subseteq A\cap L_1+A\cap L_2$ and that it has
codimension $1$. Observe, however, that the whole argument also holds
for every
$B\in\B $ such that $A\subseteq B$, and hence:
$$\limil{A\subseteq B}B\cap L_3\subseteq \limil{A\subseteq B} B\cap
L_1+\limil{A\subseteq B} B\cap L_2$$
 is of codimension $1$. From this, one concludes that $L_3\subseteq
L_1+L_2$ has codimension $1$, as desired.

Let us now compute the codimension of the other inclusion. Observe
again that the whole construction can be repeated with $B\in\B $ such
that
$B\cap L_i=(0)$ instead of $A$, and:
$$\bigwedge^k\big(V/(B\cap L_1+\dots+B\cap L_j)\big)$$
instead of $\Lambda_j$. In this case the analogous inclusion of
({\ref{eq:inclu}}) implies that:
$$(B+L_1)\cap (B+L_2)\subseteq (B+L_3)$$
is of codimension $1$. However, we can assume that $B\cap
(L_1+L_2)=(0)$, since $L_3\subseteq L_1+L_2$ has codimension $1$.
One then concludes that $L_1\cap L_2\subseteq L_3$ and that it
is of codimension $1$.

Conversely, given $L_1,L_2,L_3$ in the hypothesis, we have that
$L_1\cap L_2,L_1+L_2\in\grv $. This implies that:
$$X:=\grass\big(L_1+L_2/L_1\cap L_2\big)\,\hookrightarrow\grb(V) $$
The condition on the codimensions implies that $L_i$ lies in
$X^0:= X\cap\grv$ (for $i=1,2,3$) and that $X^0$ is a projective
line.
\end{pf}

\begin{cor}
Let $X\subset\grv$ be a projective line, and let $L_1,L_2$
be two distinct points of it. It then holds that:
$$X\,=\,\grass^1\left((L_1+L_2)/(L_1\cap L_2)\right)$$
Moreover, $L_1+L_2, L_1\cap L_2$ do not depend on the choice of
$L_1,L_2$.
\end{cor}

The same methods can easily be generalized to prove the main
result of this subsection and it is the following characterization of
finite dimensional projective spaces in (finite or infinite)
Grassmannians. Recall that $n+2$ points of a projective space define
a reference in it iff there is no  $n+1$ of them lying in a
$(n-1)$-dimensional projective space.

\begin{thm}\label{thm:proj-sp}\hfill
\begin{itemize}
\item Let  $\{L_i\}_{1\leq i\leq n+2}$ be points of $\grv$ defining a
$n$-dimensional reference (as points in $\P\Omega^*$) and assume:
$$\dim_k(L_1+\dots+L_{n+2})/(L_1\cap\dots\cap L_{n+2})=n+1$$
It then holds that:
$$\grass^k(L_1+\dots+L_{n+2})/(L_1\cap\dots\cap L_{n+2})$$
($k$ being $\dim_k L_i/(L_1\cap\dots\cap L_{n+2})$) is a
$n$-dimensional projective space contained in $\grv$.
\item If $X=\P_n\subseteq\grv$, then there exists a reference
$\{L_i\}_{1\leq i\leq n+2}$ in $X$
such that  $\dim_k(L_1+\dots+L_{n+2})/(L_1\cap\dots\cap L_{n+2})=n+1$
and:
$$X=\grass^k(L_1+\dots+L_{n+2})/(L_1\cap\dots\cap
L_{n+2})\subseteq \grv$$
where $k=1$ or $k=n$. (Note that $k$ does not depend on $\{L_i\}$ but
 on $X$ only).
\end{itemize}
\end{thm}

\subsection{Automorphisms of $\grv$}
\label{sec:aut-grfin}

In the sequel, $\aut_{\text{$k$-scheme}}(\grv)$ will be simply denote
by
$\aut(\grv)$ and similarly for $\gr^\bullet(V)$.

\begin{lem}\label{lem:auto-proy}
Let $X\subseteq\grv$ be a finite dimensional projective space and
$\phi$ be an automorphism of $\grv$. Then,
$\phi(X)$ is a finite dimensional projective space.
\end{lem}

\begin{pf}
An automorphism $\phi$ of $\grv$ induces automorphism of the
universal submodule, $\phi^*\L\iso\L$, from which one deduces
$\phi^*\detd_V\iso\detd_V$, and hence a projectivity $\phi^*$ of
$\P H^0(\grv,\detd_V)^*$. Now, the result follows easily.
\end{pf}

\begin{lem}\label{lem:aut-ext}
Fix $\phi\in\aut(\grv)$. There then exists an unique
$\bar\phi\in\aut(\gr^\bullet(V)$ with the following properties:
\begin{enumerate}
\item it is an extension of $\phi$ ($\bar\phi\vert_{\grv}=\phi$),
\item $\bar\phi$ is an inclusion-preserving or inclusion-reversing
automorphism.
\end{enumerate}
\end{lem}

\begin{pf}
Let us first define $\bar\phi(L)$ for $L\in\gr^k(V)$ ($k>0$).
Choose $L'\in\gr^{-1}(V)$ such that
$L'\subset L$, and hence
$\P(L/L')\subseteq\grv$ is a finite dimensional projective space.
Theorem {\ref{thm:proj-sp}} implies that there exist a finite
family of points $\{M_i\}$ of $\P(L/L')$, such that:
$$\P(L/L')\,=\, \P\big(\cup M_i/\cap M_i\big)$$

Using this theorem again and Lemma {\ref{lem:auto-proy}}, it follows
that:
$$\phi\big(\P(L'/L)\big)\,=\,\grass^r\big(\cup \phi(M_i)/\cap
\phi(M_i)\big)$$
where $r=1$ or $r=-1$ (codimension $1$). Observe that $r$ does not
depend on
$L$ but on the connected component of $\gr^\bullet(V)$ where $L$ lies,
that is, $r$ only depend on $k$. Now
take $L''\in\gr^{k'}(V)$ such that $L\subset L''$ (and $k< k'$),
then $\P(L/L')\subset \P(L''/L)\subseteq\grv$, and it follows that
$r$ does not depend on $k>0$. An analogous argument shows that it
does not depend on $k<0$ neither, that is, it only depends  on $\phi$.

Continuing with the above notations, we define:
$$\bar\phi(L')\,:=\, \cases
\cup \phi(M_i) & \text{ if $r=1$} \\
\cap \phi(M_i) & \text{ if $r=-1$}
\endcases$$
(observe that this definition does not depend on the choice of
$\{M_i\}$). It is now clear that $\bar\phi$ is the desired extension
which is inclusion-preserving when $r=1$ and inclusion-reversing
when $r=-1$.
\end{pf}

\begin{cor}
Let $V^+$ be finite dimensional. Then the following conditions
are equivalent:
\begin{enumerate}
\item there exists $\phi\in\aut(\grv)$ with an inclusion-reversing
extension,
\item $V$ is finite dimensional and $\dim_kV=2 \dim_kV^+$,
\item there is an irreducible hemisymmetric metric in
$V$, such that $V^+$ is m.t.i.~.
\end{enumerate}
\end{cor}

\begin{pf}
Conditions 2 and 3 are clearly equivalent. The third implies the
first since $R$, the automorphism of $\grv$ constructed in
({\ref{eq:Raut}}), extends naturally to an inclusion-reversing
automorphism of
$\gr^\bullet(V)$.

Let us prove that 1 implies 2. Let $\phi\in\aut(\grv)$ and let
$\bar\phi\in\aut(\gr^\bullet(V)) $ be its inclusion-reversing
extension. Observe that:
$$\grv=\grass^{-k}(V)$$
where $k=\dim_k(V^+)$. Since $\bar\phi$ reverses inclusions and
leaves $\grv$ invariant, it follows that:
$$\bar\phi:\grass^{-k-r}(V)\iso \grass^{-k+r}(V)\qquad \forall
r\in\Z$$
Observe that for $r=-k-1$, the scheme on the left hand side is a
projective space. By Theorem {\ref{thm:proj-sp}}, one has that
$\dim_kV=2k$.
\end{pf}

We now arrive at the classical results on the automorphism group of
Grassmannians:

\begin{cor}\label{cor:aut-fin-grass}
The group $\aut(\grv)$ is canonically isomorphic to:
\begin{itemize}
\item $\P GL(V)$ if $\dim_kV<\infty$ and $\dim_kV\neq 2 \dim_kV^+$;
\item $\P GL(V)\times \Z/2$ if $\dim_kV<\infty$ and $\dim_kV= 2
\dim_kV^+$;
\item $\P GL(V)$ if $\dim_kV=\infty$ and $\dim_kV^+<\infty$;
\item $\P GL(V^c)$ if $\dim_kV=\infty$ and $\dim_k(V/V^+)<\infty$.
\end{itemize}
\end{cor}

\begin{rem}
These results have been already proved in three remarkable papers on
the subject; namely \cite{Chow}, \cite{C} and \cite{K}. The two first
statements were proved by Chow in
\cite{Chow}, while the last two are algebraic versions of the
results of Cowen (\cite{C}) and Kaup (\cite{K}), which are given for
Hilbert spaces.

While Chow and Cowen's techniques
(\cite{Chow,C}) are based on the study of  Schubert cycles and
adjacency, our study deals with the
structure of finite dimensional projective spaces  contained in
infinite Grassmannians. Nevertheless,
such a study does imply ``adjacency type'' properties.

Recall that such a group isomorphism can not hold for
arbitrary infinite Grassmannians (see \S1 of \cite{K} for the case of
a Banach space).
\end{rem}

\begin{rem}
Recall that a collineation of a vector space over a field $k$ is a
semi-linear transformation; that is, there is an automophism of $k$
involved. If $E$ is a finite dimensional vector space, then
$\aut(\P E^*)$ consists of semi-linear transformations.
However, if we consider only $\spk$-automorphisms then the
automorphism of $k$ has to be the identity, and hence the
automorphism group of
$\P E^*$ as $\spk$-scheme is $\P Gl(E)$. Bearing this in mind,
Corollary {\ref{cor:aut-fin-grass}} can easily be generalized to the
case of $\sp(\Z)$-schemes.
\end{rem}

Our goal now is to give a similar result for the
infinite-dimensional case: $V=k((z))$. In view of the following
Proposition we can restrict our study to automorphisms with an
inclusion-pre\-ser\-ving extension, with no loss of generality.

\begin{prop}
The subset:
$$\left\{\text{ automorphisms with an inclusion-preserving extension
}\right\}$$
is a index $2$ subgroup of $\aut(\grv)$. In the case of
$\gr(k((z)))$, the quotient is generated by the automorphism $R$
of ({\ref{eq:Raut}}). (On $k((z))$ the metric:
$T_2(f,g):=\res_{z=0}f(z)g(-z)dz$ satisfies the hypothesis of
\ref{subsubsect:metric}).
\end{prop}

Recall that these two types of automorphisms correspond to
collineations and correlations in the finite dimensional
case. Nevertheless, study of the automorphism group is rather
complicated, and we shall therefore add some extra structure to the
pair $(V,\B)$. This consists of a separated linear topology on $V$
with a basis
${\frak C}$ such that:
\begin{itemize}
\item $A\sim B$  for all $A,B\in{\frak C}$,
\item $A\in\gr^\bullet(V)$ for all $A\in{\frak C}$.
\end{itemize}
In this setting $\;\check{}\;$ will denote the completion w.r.t.
the topology induced by ${\frak C}$.

Now, we define:
$$\begin{gathered}
\gl({\frak C})\,:=\,
\left\{\begin{gathered}
\phi\in\glv(k)\text{ bicontinuous w.r.t. }{\frak C} \\
\text{ such that } \phi(\grv)=\grv \end{gathered}\right\}
\\
{\mathcal A}({\frak C})\,:=\,
\left\{\begin{gathered}
\phi\in\aut(\grv )\text{ bicontinuous w.r.t. $\frak C$,}\\
\text{ with an inclusion-preserving extension }\end{gathered}
\right\}
\end{gathered}$$
where bicontinuous w.r.t. $\frak C$ means that both
$\phi(A)$ and $\phi^{-1}(A)\in{\frak C}$ contain an element of
${\frak C}$ (for all $A\in{\frak C}$).

Since $k^*$ is the kernel of the morphism $\glv(k)\to\aut(\grb(V))$,
one has the exactness of the following sequence:
\beq
0\to k^*\to \gl({\frak C})\to {\mathcal A}({\frak C})
\label{eq:gl-aut}\end{equation}

The following result is to be interpreted as the natural
generalization of Corollary {\ref{cor:aut-fin-grass}} to the infinite
dimensional case.

\begin{thm}\label{thm:pgl-aut-pgl}
There exist injective morphisms of groups:
$$\P\gl({\frak C})\,\hookrightarrow\,
{\mathcal A}({\frak C})\,\hookrightarrow
\P\gl({\frak C})(\check V)$$
such that the composite maps a continuous automorphism of $V$ to
that canonically induced on $\check V$.
\end{thm}

\begin{pf}
Since $\bar\phi$ is inclusion-preserving, we have the
following commutative diagram:
\beq
\begin{array}{ccccc}
\gr(V/L) & \hookrightarrow & \gr(V/L') & \hookrightarrow & \grv
\\ \downarrow & & \downarrow & & \downarrow  \\
\gr(V/\phi(L)) & \hookrightarrow & \gr(V/\phi(L'))
& \hookrightarrow & \grv \end{array}
\label{eq:comm-dia}\end{equation}
where $L'\subseteq L$ and $L,L'\in{\frak C}$. Applying Corollary
{\ref{cor:aut-fin-grass}}, one obtains a family of linear
applications (up to a scalar):
$$T_L:V/L @>{\,\sim\,}>>V/\phi(L)$$
Bearing in mind the inclusion  $L'\hookrightarrow L$
and ({\ref{eq:comm-dia}}), the commutativity of the following diagram
becomes apparent:
$$\CD V/L @>T_L>> V/\phi(L)\\
@VVV @VVV \\
V/L' @>T_{L'}>> V/\phi(L')\endCD$$
or, what amounts to the same, $\{T_L\}_{L\in{\frak C}}$ is a
morphism of inverse systems.  Taking inverse limits one obtains
$T_\phi:\check V\iso \check V$, which is bicontinuous w.r.t.
${\frak C}$ by the very construction.

Let us show that if $T_\phi=\lambda\cdot Id$ ($\lambda\in k^*$), then
$\phi=Id$. Note first that $T_\phi=Id$ implies that $\phi$ induces the
identity on $\gr(V/L)$ for all $L\in{\frak C}$, and hence
the commutative diagram of affine schemes:
$$\CD F_A @>{\phi}>> F_A \\ @A{j_L}AA @AA{j_L}A \\
j_L^{-1}(F_A) @>Id>> j_L^{-1}(F_A) \endCD$$
where $j_L:\gr(V/L)\to\grv$. Taking inverse limits in the inverse
system of morphisms between the  associated rings, it is easy to
conclude that $\phi\vert_{F_A}=Id\vert_{F_A}$ and hence $\phi=Id$,
since there is no non-constant function on $F_A$ vanishing on
$F_A\cap \gr(V/L)$ for all $L\in{\frak C}$.
\end{pf}

\begin{exam}
Assume that $\dim(V/V^+)$ is finite and let ${\frak C}$ be the set of
all rational points of $\grv $. It then holds that $\P\glv\iso
\P\gl({\frak C})$ and ${\mathcal A}({\frak C})\iso
\aut(\grv )$. Summing up:
$$\P\glv\,\iso\,\aut(\grv )$$
\end{exam}

\begin{exam}
It is clear from Theorem {\ref{thm:pgl-aut-pgl}} that, in general, it
is not possible to represent elements of the linear groups above as
matrices. Nevertheless, let us study the case of $V=k((z))$.

Let ${\frak C}$ be the set of subspaces $\{z^n\cdot
k[z^{-1}]\}_{n\in\Z}$. One checks that $\check V=k[[z^{-1},z]]$, and
hence to $g\in\gl(\check V)$ there is an
associated $\Z\times\Z$-matrix $(g_{ij})_{i,j\in\Z}$ ($g_{ij}$
being the coefficient of $z^j$ in $g(z^i)$). But one recovers $g$
from the matrix if and only if $g$ is continuous w.r.t. the
$z$-adic topology.

Moreover, the image of:
$$\gl({\frak C})\longrightarrow
M_{\Z\times\Z}(k)$$
consists of matrices $(g_{ij})_{i,j\in\Z}$ such that there exists
$n(j):\Z\to\Z$ satisfying: $g_{ij}=0$ for $i<n(j)$, and $g_{n(j)j}\in
k^*$.   Further, it is clear
that the matrices of the group $\gl_{res}$ of Segal-Wilson belong to
the image.
\end{exam}


\end{document}